\documentclass[12pt]{article}
\usepackage{amsmath,amssymb,amsfonts}

\newcommand{\R}{\mathbb{{R}}}
\newcommand{\Z}{\mathbb{{Z}}}

\let\vp=\varphi
\let\st=\stackrel
\let\l=\longrightarrow
\newcommand{\para}{\paragraph}
\begin{document}
\begin{center}
\vspace{1cm}
{\bf\large On the Poincar\'{e} Index of Isolated Invariant Sets}

\vspace*{5mm}

 M.R. Razvan  and  M. Fotouhi Firoozabad\\

{\it \small Institute for Studies in Theoretical Physics and Mathematics\\
P.O.Box: $19395-5746$, Tehran, IRAN\\  E-Mail:
razvan@karun.ipm.ac.ir\\
Department of Mathematics, Sharif University of
Technology\\  P.O.Box: $11365-9415$, Tehran, IRAN\\
E-Mail: fotuhifi@math.sharif.ac.ir\\
E-Mail: fotouhi@karun.ipm.ac.ir\\}
\end{center}

\begin{abstract}
In this paper, we use Conley index theory to examine the
Poincar\'{e} index of an isolated invariant set. We obtain some
limiting conditions on a critical point of a planar vector field
to be an isolated invariant set. As a result we show the existence
of infinitely many homoclinic orbits for a critical point with
the Poincar\' {e} index greater than one.
\end{abstract}


\section{Introduction}

The Conley index has proved to be a useful tool in the
investigation of qualitative properties of dynamical systems. It
has generalized Morse theory for an isolated invariant set of a
continuous flow on a locally compact metric space\cite{C,CZ2}.
For this reason, Conley index is known as a generalization of
Morse theory.  In \cite{CZ1} and \cite{CZ2}, Conley and Zehender
used this index to show the existence of periodic solutions for
Hamiltonian systems. This was a landmark in the proof of the
Arnold conjecture on the existence of periodic orbits of
Hamiltonian systems on symplectic manifolds. Conley index theory
has also some applications in the existence of solutions for a
class of differential equations. (See \cite{Ry,Sm,MM} and
references therein.)

In Conley index theory, we deal with pairs of closed sets called
index pair for an isolated invariant set $I$. The homotopy type
of these index pairs is independent of the index pair chosen,
which is called the Conley index of $I$ and denoted by $h(I)$.
This paper concerns the relation between Conley index theory and
the Poincar\'e-Hopf theorem \cite{Mc1}. We define the Poincar\'e
index of an isolated invariant set $I$ to be $\chi (h(I))$. This
definition coincides with the classical Poincar\'e index when the
invariant set is a single point. We use some topological
properties of the Conley index to obtain restrictions on the
Poincar\'e index of isolated invariant sets in dimension two. It
is well-known that on a two-dimension manifold M, the Poincar\'e
index of an isolated critical point of a gradient vector field is
not greater than one. Here we show that a critical point $x$ with
$ind(x)>1$ cannot be an isolated invariant set. This concludes the
existence of infinitely many homoclinic orbits for such a
critical point.

We first present some basic results from Conley index theory. Then
we define the concept of continuation and provide a new proof for
the results of \cite{Mc1} based on continuation to gradient
\cite{Re1,Re2}. Finally we apply these results to show the
existence of infinitely many homoclinic orbits in dimension two.

\section{Conley Index}

Let $\vp^t$ be a $C^1$-flow on a smooth manifold $M$. A subset
$I\subset M$ is called an {\it isolated invariant} set if it  is
the maximal invariant set in some {\bf compact} neighborhood of
itself. Such neighborhood is called an {\it isolating
neighborhood}.

\para{Definition.} A closed pair $(N,L)$ is called an {\it index pair}
for $I$ if
\begin{enumerate}
\item $\overline{N-L}$ is an isolating neighborhood for I.
\item $L$ is positively invariant relative to $N$, i.e., if $x\in L$,
$t\geq 0$, $\vp^{[0,t]}(x)\subset N$, then $\vp^{[0,t]}(x)\subset L$.
\item $L$ is the exit set of $N$, i.e., if $x\in N$, $t\in\R^+$ and
$\vp^t(x)\notin N$, then there is a $t'\in [0,t]$ such that
$\vp^{t'}(x)\in L$.
\end{enumerate}

In \cite{C,CZ2,Sa,Sm} it has been shown that every isolated
invariant set $I$ admits an index pair $(N,L)$ and the homotopy
type of $(N/L,[L])$ is independent of the index pair chosen. We
denote the homotopy type of $(N/L,[L])$ by $h(I)$ and call it the
{\it Conley index} of $I$. The {\it homology Conley index} of $I$
is defined by $CH_*(I)=H_*(N/L,[L])$.

\para{Example 2.1.} Let $x\in M$ be a nondegenerate
critical point for $f:M\st{C^2}{\l}\R$. Then $\{x\}$ is an
isolated invariant set for $-\nabla f$ and by Morse Lemma
\cite{M}, it is easy to show that $h(\{p\})$ is a pointed
k-sphere where $k$ is the number of positive eigenvalues of
Hessian matrix $f$ at $p$. Therefore the Conley index can be
considered as a generalization of Morse index.

It is not true that $H_*(N,L)\cong H_*(N/L,[L])$ for every index
pair $(N,L)$. In \cite{Sa}, Salamon introduced a class of index
pairs for which the above isomorphism holds.
\para{Definition.} An index pair $(N,L)$ is called {\it regular} if the
exit time map $$\tau_+:N \l [0,+\infty] , \ \ \tau_+(x)=\left\{
\begin{array}{ll}
sup\{t|\varphi^{[0,t]} (x) \subset N-L\} & \text{if} \ x\in N-L,\\
0 & \text{if}\ x\in L,
\end{array} \right. $$
is continuous. ( See \cite{Sa} for more details about regular
index pairs.) For every regular index pair $(N,L)$, we define the
induced semi-flow on $N$ by $$ \varphi ^t _{\natural} :N\times \R
^+ \l N,\ \ \varphi^t _{\natural}(x)=\varphi ^{min\{t,\tau_+
(x)\}} (x)$$

\para{Proposition 2.2.} If $(N,L)$ be a regular index pair for a
continuous flow $\varphi ^t$, then $L$ is a neighborhood
deformation retract in $N$. In particular, the natural map $\pi :
N \l N/L$ induces an isomorphism $H_*(N,L)\cong H_*(N/L,[L])$.

\noindent{\bf Proof.} Consider the induced semi-flow
$\varphi_{\natural}$ on $N$ and the neighborhood $U:=\tau _+
^{-1} [0,1]$ of $L$. Now $\varphi_{\natural}|_{U\times [0,1]}$
gives the desired deformation retraction. $\square$

In \cite{RS}, Robbin and Salamon proved that every isolated
invariant set admits a regular index pair which is stable under
perturbation. They first showed the existence of a smooth
Liapunov function on a neighborhood of the isolated invariant set.

\para{Theorem 2.3.} Let $N$ be an isolating neighborhood of
$I$. Then there is a neighborhood $U$ of $N$ and a smooth function
$f:U{\l}\R$ satisfying
\begin{itemize}
\item [(i)] $f(x)=0$ for all $x\in I$.
\item [(ii)] $\frac{d}{dt}|_{t=0}f(\vp^t(x))<0$ for all $x\in N-I$.
($f$ decreases along orbits in $U-I$.)
\end{itemize}

Then they used this Liapunov function to construct a triple
$(N,L^-,L^+)$, such that $(N,L^+)$ is a regular index pair for
$I$ with respect to the forward flow and $(N,L^-)$ is a regular
index pair for $I$ with respect to the reverse flow. Furthermore
$L^+$ and $L^-$ can be chosen to be (n-1)-manifolds with
boundary, so that $N$ is a manifold with corners with those
corners contained in $L^-\cap L^+$ and $N=L^-\cup L^+$. We call
such a triple $(N,L^-,L^+)$ as a {\it regular index triple} for
$I$ in $M$. The Conley indices of $I$ related by the forward and
reverse flow are represented by $h^+(I)$ and $h^-(I)$. If $M$ is
 orientable in a neighborhood of $I$, the indices for the forward and reverse
flows are related by Poincar\'{e}-Lefschetz duality isomorphism
$H_*(N,L^+)\simeq H^{m-*}(N,L^-)$ where $m=dim M$. (See
\cite{Mc2,MS,Sp}.) If we consider the homology with coefficients
in $ \Z_2$, the Poincar\'e-Lefschetz duality is valid without the
assumption of orientability.

\para{Definition.} $A\subset M$ is called an {\it attractor} set if
it is the $\omega$-limit set of a compact neighborhood of itself.
A {\it repeller} set is an attractor set for the reverse flow.

\para{Proposition 2.4.} $I$ is an attractor set for $\vp^t$ if
and only if there is an index pair $(N,L)$ for $I$ which
$L=\varnothing$.

\noindent{\bf Proof.} Let $I$ be an attractor and $V$ be a
neighborhood of $I$ such that $\omega (V)=I$. Then there is a
$T>0$ such that $\vp^{[T,\infty)}(V)\subset int(V)$. If we set
$N:=\bigcap_{0\leq s\leq T}\vp^s(V)$, then $(N,\varnothing)$ is
an index pair for I. Now assume that $(N,\varnothing)$ is an index
pair for $I$. According to the property (3) of the definition of
index pair, we imply that $N$ is positively invariant, hence
$\omega(N)\subset N$. Since  $N$ is an isolating neighborhood for
$I$, it follows that $\omega(N)\subseteq I$. Since $I$ is an
invariant set, we conclude that $\omega(N)=I$. $\square$

\para{Theorem 2.5.} Suppose that $I\subset M$ is a connected
isolated invariant set.
\begin{itemize}
\item [(i)] If $I$ is not an attractor, then $H_0(h^+(I))=0$.
\item [(ii)] If  $I$ is
not a repeller, then $H_m(h^+(I);\Z_2)=0$. Moreover if $M$ is
orientable in a neighborhood of $I$, then $H_m(h^+(I))=0$.
($m=\dim M$)
\end{itemize}

\noindent{\bf Proof.} Consider a regular index triple
$(N,L^+,L^-)$ for $I$. We may assume that $N$ is connected
(otherwise replace $N$ by the connected component  of $N$ that
contains $I$). Since $I$ is not an attractor set,
$L^+\neq\varnothing$ by Proposition 2.4. Thus
\[
H_0(h^+(I))=H_0(N,L^+)=0
\]
Similarly we have $L^-\neq\varnothing$ and $H^0(h^-(I))=0$. Now by
Poincar\'{e}-Lefschetz duality $H_m(h^+(I))\simeq H^0(h^-(I))=0$.
$\square$

\para{Definition.} A {\it Morse decomposition} of $I$ is a finite collection
$\{M_i\}_{i=1}^n$ of disjoint, nonempty isolated invariant
subsets of $I$ such that for each $x\in I-(\bigcup_{i=1}^nM_i)$,
there are $1\leq i<j\leq n$ such that $\alpha(x)\in M_j$ and
$\omega(x)\in M_i$.

\para {Example 2.6.} Consider the gradient flow of a smooth
function $f$ on a compact manifold $M$. Suppose that
$\{x_1,\cdots,x_n\}$ are critical points of $f$ with $f(x_i)\leq
f(x_j)$ for $i<j$. Then $\{x_1,\cdots,x_n\}$ is a Morse
decomposition for $M$.

\para{Theorem 2.7.} Let $I$ be an isolated invariant set with a
Morse decomposition $\{M_i\}_{i=1}^n$. Define the Poincar\'{e}
polynomial $p(t,X):=\sum_{k=0}^\infty b_k(X)t^k$ that
$b_k(X)=\dim(H_k(X))$ is k-th Betti number of $X$. Then there is
a polynomial $Q$ with nonnegative coefficients such that
\[
\sum_{i=1}^np(t,h(M_i))=p(t,h(I))+(1+t)Q(t)
\]

The above theorem is known as the generalized Morse inequalities
\cite{CZ2,Sm}. If we apply this result to the gradient flow of a
Morse function with the Morse decomposition described in Example
2.6., we obtain the classical Morse inequalities \cite{M}.

\section{Continuation}

A parametrized flow on $M$ is a collection of flows
$\{\vp_t^\lambda\;|\;\lambda\in I\}$ indexed by $I=[0,1]$ such
that $\Phi_t(x,\lambda)=(\vp_t^\lambda(x),\lambda)$ is a
$C^1$-flow on $M\times I$. we say $S^0$, an invariant set for
$\vp_t^0$, and $S^1$, an invariant set for $\vp_t^1$, are related
by {\it continuation} if there is an isolated invariant set
$\Sigma\subset M\times I$ for $\Phi_t$ such that $S^0=\Sigma\cap
M\times\{0\}$ and $S^1=\Sigma\cap M\times\{1\}$. The reason that
continuations are interesting in Conley index theory is the
following theorem \cite{C,Sm}.

\para{Theorem 3.1.} If $S^0$ and $S^1$ are isolated invariant sets
related by continuation, then $h(S^0)=h(S^1)$.

f we have an isolated invariant set $S$ and make a small
perturbation of the flow, then the new flow will have an isolated
invariant set $S'$ near $S$, and by the above theorem
$h(S)=h(S')$. Therefore Conley index is invariant under
perturbation. In \cite{Re1,Re2} Reineck has shown that every
isolated invariant set can be continued to an isolated invariant
set in a Morse-Smale gradient flow. Since Morse-Smale flows are
easy to deal with, we first prove our result for a Morse-Smale
flow and then extend it to the general case by using the
invariance of Conley index  under continuation.

\para{Theorem 3.2.} Let $X$ be a smooth vector field on a
Riemannian manifold $M$ and let $I$ be an isolated invariant set
in $\vp_t$, the flow generated by $X$, with isolating
neighborhood $N$. Then $I$ can be continued to an isolated
invariant set in a Morse-Smale gradient flow without changing the
vector field on $M-N$.

\para{Definition.}  We define the Poincar\'{e} index of an
isolated invariant set $I$ to be the Euler characteristic of the
Conley index of $I$, i.e. $ind_p(I):=\chi(h(I))$.

Suppose that the flow $\vp^t$ is associated with a vector field
$X$ on $M$. If $\{x\}$ is a critical point of $X$ and an isolated
invariant set for $\vp^t$, then $ ind_p(x)$ coincides with the
classical definition of Poincar\'{e} index of $x$ (up to a sign).
This is a special case of the results of \cite{Mc1} in which
McCord developed the Poincar\'e-Hopf theorem and showed that
$ind_p(I)=(-1)^m\sum ind(x)$, where the sum is taken over all
critical points in $I$, $ind(x)$ is the Poincar\'{e} index of $x$
relative to vector field $X$ and $m=dimM$. We provide another
proof for this result based on Reneik's  continuation to gradient.

\para{Theorem 3.3.} Let $I$ be an isolated invariant set. Then
\[
ind_p(I)=(-1)^m\sum_{x\in I}ind(x).
\]
In particular, if $ind_p(I)\neq 0$, then there exists a critical
point in $I$.

\noindent{\bf Proof.} By Theorem 3.2., $I$ can be continued to an
isolated invariant set $J$ in a Morse-Smale gradient flow
$-\nabla f$ without changing the vector field on $M-N$ where $N$
is an isolating neighborhood for $I$. Thus $h(I)=h(J) $ by
Theorem 3.1. If $\{y_i\}_{i=1}^n$ is the set of critical points
of $f$ in $N$, then $\sum_{x\in I} ind(x)=\sum_{i=1}^nind(y_i)$.
If the Morse index of $y_i$ is $k$, then $ind(y_i)=(-1)^{m+k}$
and by Example 2.1, $h(\{y_i\})=\Sigma^k$. Let $\mu_k$ denotes
the number of critical points of Morse index $k$. Thus
$\sum_{x\in I} ind(x)=\sum_k(-1)^{m+k}\mu_k$. By Example 2.6.,
$\{y_i\}_{i=1}^n$ is a Morse decomposition for $J$ with respect
to $-\nabla f$. According to the generalized Morse inequalities
(Theorem 2.11), we have
\[
\sum_{i=1}^np(-1,h(y_i))=p(-1,h(J))=p(-1,h(I)).
\]
Since $h(y_i)=\Sigma^k$, we get $p(t,h(y_i))=t^k$ and
$\sum_k(-1)^k\mu_k=p(-1,h(I))$. Therefore
\[
(-1)^m\sum_{x\in I} ind(x)=p(-1,h(I))
=\sum_{k=0}^{\infty}(-1)^kb_k(h(I))=\chi(h(I))=ind_p(I).
\hspace{0.5cm}\square
\]

\para{ Proposition 3.4.} Suppose that $I\subset M$ is an NDR
(Neighborhood Deformation Retract) isolated invariant set.
\begin{itemize}
\item [(i)] If $I$ is an attractor, then $ind_p(I)=\chi(I)$.
\item [(ii)] If  $I$ is a repeller, then $ind_p(I)=(-1)^m\chi(I)$.
  ($m=\dim M$)
\end{itemize}

\noindent{\bf Proof.} When $I$ is an attractor, there is an index
pair $(N,\varnothing)$ for $I$ by Proposition 2.4. Since $I$ is an
NDR, there exists a neighborhood $U\subset N$ such that $I$ is
deformation retract of $U$. By the definition of index pair, $N$
is positively invariant and $\omega(N)=I$. So there is a $T>0$
such that $\vp^T(N)\subset U$. Therefore $N$ can be deformed to
$I$ and $H_i(N,\varnothing)=H_i(I)$ for every $i$. Thus
$\chi(h(I))=\chi(I)$ which proves (i). Notice that for a finite
CW-complex, the Euler characteristic does not depend on the
coefficients field. Since $I$ is assume to be an NDR, it has the
homotopy type of a finite CW-complex. If we consider the homology
with coefficients in $\Z_2$, we obtain $\chi(h^+(I))=(-1)^m
\chi(h^-(I))$ by the duality theorem. So if $I$ is an NDR
repeller, then $ind_p(I)=(-1)^m\chi(I)$. $\square$

\section{Applications}

In this section, we consider a smooth vector field on a {\bf
surface} $M$ with an isolated critical point $x$. We desire to
show that if $ind(x)>1$, then $x$ is accumulated by infinitely
many homoclinic orbits.
\para{Lemma 4.1.} Let $I\subset M$ be a connected NDR
isolated invariant set such that $ind_p(I)>0$. Then $I$ is either
an attractor or a repeller and $ind_p(I)=\chi(I)$.

\noindent{\bf Proof.} Suppose that  $I$ is neither an attractor
nor a repeller.  By Theorem 2.5., $H_2(h(I);\Z_2)\simeq
H_0(h(I);\Z_2)=0$. Now we conclude that
\[\begin{array}{lll}
ind_p(I)&=&\chi(h(I))\vspace{2mm}\\
&=&rank(H_2(h(I);\Z_2))-rank(H_1(h(I);\Z_2))+rank(H_0(h(I);\Z_2))\vspace{2mm}\\
&=&-rank(H_1(h(I);\Z_2))\leq 0.
\end{array}\]
Since $m=2$, the proof is complete by Proposition 3.4. $\square$

\para{Theorem 4.2.}  Let $x$ be a critical point for a vector field
on a surface $M$. If $ind(x)>1$, then there exists a homoclinic
orbit in any neighborhood of $x$.

\noindent{\bf Proof.} We first show that $\{x\}$ cannot be  an
isolated invariant set. Suppose the contrary, then according to
Theorem 3.3. and the above lemma, $ind(x)=\chi(\{x\})=1$ which is
a contradiction. Consider a closed neighborhood $V$ of $x$ with no
critical points rather than $x$. Let $I(V)$ be the maximal
invariant set in $V$. The above argument says that there is a
point $y\neq x$ in $I(V)$, hence $\omega(y)\subset I(V)$. Notice
that there cannot be any cycle in $V$. To see this, suppose that
is $\gamma$ is a cycle in $V$. Then the Poincar\'{e} index of
$\gamma$ must be one \cite{Pe}, thus there exists a critical
point inside $\gamma$.  Since the only critical point in $V$ is
$x$, we get $ind(x)=1$ which is a contradiction. Now according to
the Poincar\'{e}-Bendixon theorem \cite{Pe,Pd}, $\omega(y)$ and
$\alpha(y)$ are critical points or homoclinic orbits. If neither
of $\omega(y)$ and $\alpha(y)$ are homoclinic orbits, then
$\omega(y)=\alpha(y)=x$. So there exists a homoclinic orbit in
$V$. $\square$

\para{Proposition 4.3} Let $\gamma$ be a homoclinic orbit with no
critical point inside of it. Then all the orbits inside $\gamma$
are homoclinic.

\noindent{\bf Proof.} Let $x:=\omega(\gamma)=\alpha(x)$ and
$\Omega$ be the region surrounded by $\gamma$. Similar to the
above argument, there is no cycles in $\Omega$ and the limit sets
of any orbit in $\Omega$ are either $\{x\}$ or homoclinic orbits.
Since $x$ is the only critical point in $\Omega$, it belongs to
all limit sets. On the other hand, it is known that if one of the
limit sets is not a critical point, then the limit sets are
disjoint \cite{Pe,Pd}. Therefore the limit sets of any orbit in
$\Omega$ must be $\{x\}$. $\square$

\para{Remark 4.4.} It is well-known that if $x$ is a critical
point of a gradient vector field, then $ind(x)\leq 1$. ( See
\cite{Pr} for another proof.) The above theorem clearly shows why
this result is true.

\para{Acknowledgment.} The authors would like to thank Institute
for studies in theoretical Physics and Mathematics, IPM , for
supporting this research.

\end{document}